\documentclass[11pt, leqno]{article}
\usepackage{amsmath, amsfonts}
\begin{document}

\author{Ajai Choudhry and Jaros\l aw Wr\'oblewski }
\title{An Ancient  Diophantine Equation \\ with applications to \\Numerical Curios and Geometric Series}
\date{}
\maketitle

\begin{abstract} In this paper we examine the diophantine equation $x^k-y^k=x-y$ where $k$ is a positive  integer $\geq 2$, and consider its applications. While the complete solution of the equation $x^k-y^k=x-y$ in positive rational numbers is already known when $k=2$ or $3$, till now only one numerical solution of the equation in positive rational numbers has been published when $k=4$, and  no nontrivial solution is known when $k \geq 5$. We describe a method of generating infinitely many positive rational solutions of the equation when $k=4$. We use the positive rational solutions of the equation with $k=2,\, 3$ or 4 
to produce numerical curios involving square roots, cube roots and fourth roots, and as another application of these solutions, we  show how to construct  examples of geometric series with  an interesting property.
\end{abstract}

\noindent {\bf Keywords:} cube roots; fourth roots; cubes; biquadrates;  geometric series; cubic diophantine equation. 

\noindent Mathematics Subject Classification: 11D25.

\setcounter{section}{0}

\section{Introduction}

\hspace{0.25in}It seems strange at first sight but it is nevertheless true that
\begin{multline}
\sqrt[\scriptstyle4]{\frac{62304353849776801}{1423276677734560000}}+{\frac{5497}{17270}}\\[0.2in]
=\sqrt[\scriptstyle4]{\frac{62304353849776801}{1423276677734560000}+{\frac{5497}{17270}}} 
\ \phantom{.} 
. \quad \quad \quad \quad \quad \quad \label{ex1}
\end{multline}
This is not an isolated curiosity. There exist infinitely many similar identities involving square roots, cube roots and fourth roots. Surprisingly, such identities may be obtained from positive rational solutions of the diophantine equation,
\begin{equation}
a^k-b^k=a-b,
\label{eqk}
\end{equation}
where $k$ is a given positive integer $ \geq  2$. 

When $k=2$, Eq.~\eqref{eqk} reduces to the trivial equation $a+b=1$. When $k=3$, Eq.~\eqref{eqk} dates back to Diophantus who gave a numerical solution in his Arithmetica \cite[p. 173]{He}. The complete solution of  \eqref{eqk} with $k=3$  is readily obtained and is already known. When $k=4$, a single numerical solution of Eq.~\eqref{eqk} in positive rational numbers was obtained by Fermat \cite[pp. 319-20]{He}. Till now, no other solutions of Eq.~\eqref{eqk} in positive rational numbers have been published when $k=4$. Further, no nontrivial solutions of \eqref{eqk} are known when $k \geq 5$. 

In Section 2, we examine Eq.~\eqref{eqk} and, in particular, we  describe a method of generating infinitely many positive rational solutions of the equation when $k=4$. In Section 3, we describe how identities of the type \eqref{ex1} may be derived using the solutions of \eqref{eqk}. Finally, in Section 4, we  give an application of the solutions of the equation \eqref{eqk} to the construction of  convergent geometric series 
$\sum_{n=1}^\infty a_n$ with positive rational terms and having the striking property that 
\begin{equation}
\sum_{n=1}^\infty a_n=
\sum_{n=1}^\infty a_n^k 
,\label{series}
\end{equation}

\section{The diophantine equation $a^k-b^k=a-b$}

\subsection{} For any positive integer value of $k$, the diophantine equation \eqref{eqk} always has the trivial solutions $(a,b)=(1,\,0)$, $(a,b)=(0,\,1)$ and $(a,b)=(m,\,m)$ where $m$ is an arbitrary rational number. Any solution other than these trivial solutions will be called a nontrivial solution. 

For a nontrivial solution, $a \neq b$, and so we may  cancel out the factor $a-b$ from both sides of \eqref{eqk} .  As already noted, when $k=2$, Eq.~\eqref{eqk} reduces to the trivial equation $a+b=1$ that has the obvious solution $(a,b)=(m,\,1-m)$ where $m$ is an arbitrary parameter. This  yields nontrivial solutions involving negative integers as well as solutions in positive rational numbers  when $m$ is rational and $ 0 < m < 1$. When $k > 2$, it is readily seen that Eq.~\eqref{eqk} has no nontrivial solutions in integers.

When $k=3$, Eq.~\eqref{eqk} reduces, on cancelling out the factor $a-b$, to the quadratic equation $a^2+ab+b^2=1$ whose complete solution is readily obtained, and is given by,
\begin{equation}
a=\frac{\displaystyle m^2-n^2}{m^2+mn+n^2},\quad b=\frac{\displaystyle 2mn+n^2}{\displaystyle m^2+mn+n^2}, \label{soleq2}
\end{equation}
where $m$ and $n$ are arbitrary integer parameters. When $m=3,n=1$, we get the numerical solution 
\begin{equation}
(a,\,b)=(8/13,\,7/13), \label{solDiophantus}
\end{equation}
 which is the solution given by Diophantus to which a reference was made in the Introduction.

\subsection{} We now consider the equation,
\begin{equation}
a^4-b^4=a-b \label{eqk4}
\end{equation}
obtained by taking $k=4$ in \eqref{eqk}.

The  solution of \eqref{eqk4} found by Fermat is given below:
\begin{equation}
a=26793/34540, \quad  b=15799/34540. \label{sol1eqk4}
\end{equation}
The same numerical solution of Eq.~\eqref{eqk4}, found by three different methods, is given in \cite[Part I, pp. 48-51]{Le}. Till date no other solution of \eqref{eqk4} in positive rational numbers has been published.  There is a brief comment on Eq.~\eqref{eqk4} by Scott and Styer \cite{SS}  who, referring to Skolem's book \cite{Sk}, have  stated that this equation has infinitely many solutions in positive rational numbers. 

Eq.~\eqref{eqk4} reduces, on cancelling out the factor $a-b$ to the cubic equation,
\begin{equation}
a^3+a^2b+ab^2+b^3=1, \label{cubic1}
\end{equation}
which represents an elliptic curve. On applying the birational transformation given by
\begin{equation}
a =(y+2)/(2x), \quad b = -(y-2)/(2x), \label{birat1}
\end{equation}
and
\begin{equation}
x = 2/(a+b), \quad y = 2(a-b)/(a+b), \label{birat2}
\end{equation}
the cubic equation \eqref{cubic1} reduces to the following Weierstrass  model of an elliptic curve:
\begin{equation}
y^2=x^3-4. \label{cubic2}
\end{equation}

Equations \eqref{birat1} and \eqref{birat2} establish a one-one correspondence between the rational solutions of Eq.~\eqref{cubic1} and the rational points on the elliptic curve \eqref{cubic2}. It is readily verified from Cremona's well-known tables on elliptic curves that Eq.~\eqref{cubic2} defines an elliptic curve over ${\mathbb Q}$ of rank 1 and a basis for the Mordell-Weil group is given by the point $P$ with co-ordinates $(2,2)$. 

Fermat's solution \eqref{sol1eqk4} corresponds to the point $A$ on the elliptic curve \eqref{cubic2} given by $(X,\,Y)=(785/484, 5497/10648)$. Since the curve \eqref{cubic2} has positive rank, there are infinitely many rational points on the elliptic curve \eqref{cubic2} and it  follows from   a theorem of Poincare and Hurwitz \cite[Satz 11, p. 78]{Sk} that  there are infinitely many rational points in the neighbourhood of the point $A$, and these points correspond to infinitely many solutions in positive rational numbers of Eq.~\eqref{cubic1} as stated by Scott and Styer. 

The aforementioned theorem of Poincare and Hurwitz does not give any method of generating the infinitely many rational points in the neighbourhood of the point $A$, and hence is of no help in finding solutions of Eq.~\eqref{eqk4} in positive rational numbers.  We give below a method of finding infinitely many solutions of Eq.~\eqref{eqk4} in positive rational numbers. 

It follows from Eqs.~\eqref{cubic1}, \eqref{birat1} and \eqref{birat2} that the points $(1,\,0)$ and $(0,1)$ on the curve \eqref{cubic1} correspond to the points $(2,\,2)$ and $(2,\,-2)$ respectively on the curve \eqref{cubic2} and there is a one-one correspondence between points on the curve \eqref{cubic1} with positive rational coordinates and  rational points on the curve \eqref{cubic2} with $2^{2/3} < x <  2$. We will first generate an infinite sequence of rational points on the curve \eqref{cubic2} with $2^{2/3} < x <  2$, and then use the relations \eqref{birat1} to find infinitely many positive rational solutions of \eqref{cubic1}.

Let $P_1$ be a known rational point on the curve \eqref{cubic2} with co-ordinates $(x_1,\,y_1)$ such that $2^{2/3} < x_1 <  2$ and $y_1 > 0.$ Using the group law, we now compute the point $Q=2P_1-P$ on the curve \eqref{cubic2}. The abscissa $x_2$ of the  point $Q$ is given by,
\begin{equation}
\begin{aligned}
x_2&=2\{x_1^8+8x_1^7-64x_1^6+64x_1^5+224x_1^4+512x_1^3+1024x_1^2\\
& \quad \quad \quad -1024x_1-1024+(4x_1^6-320x_1^3-512)y_1\}\\ 
& \quad \quad  \quad \quad \times (x_1^4-8x_1^3+32x_1+32)^{-2}.
\end{aligned}
\end{equation}
We will now show that $2^{2/3} <  x_2 <  2$. The first part of this inequality is readily established since $Q$ is a rational point that  lies on the curve \eqref{cubic2}, and so it follows that  $2^{2/3} <  x_2$. To prove the second part, we note that,
\begin{equation}
\begin{aligned}
2-x_2&=8\{(-x_1^6+80x_1^3+128)y_1-2(3x_1^4-16x_1^3+96x_1+64)\\
& \quad \quad \quad \times (x_1^3-4)\}(x_1^4-8x_1^3+32x_1+32)^{-2}. \label{diff1}
\end{aligned}
\end{equation}
and  since $2^{2/3} < x_1 <  2$ and $y_1 > 0,$ therefore $(-x_1^6+80x_1^3+128)y_1 > 0$ and also $2(3x_1^4-16x_1^3+96x_1+64)(x_1^3-4) > 0$. Further, since $y_1^2=x_1^3-4$, we have, 
\begin{multline}
\{(-x_1^6+80x_1^3+128)y_1\}^2-\{2(3x_1^4-16x_1^3+96x_1+64)(x_1^3-4)\}^2\\
=(x_1-2)(x_1^3-4)(x_1^3-18x_1^2-36x_1-40)(x_1^4-8x_1^3+32x_1+32)^2. \label{diff2}
\end{multline}
Since $2^{2/3} < x_1 <  2$, it is readily seen that the quantity on the right-hand side of \eqref{diff2} is positive. It follows that 
\begin{equation}
(-x_1^6+80x_1^3+128)y_1 > 2(3x_1^4-16x_1^3+96x_1+64)(x_1^3-4),
\end{equation}
and therefore, from \eqref{diff1}, we get $2-x_2 > 0$. We have thus proved that $2^{2/3} <  x_2 <  2$. 

The two points $Q$ and $-Q$ on the curve \eqref{cubic2} have the same positive abscissa $x_2$ and one of them necessarily has positive ordinate as well. We take $P_2$ to be either $Q$ or $-Q$ such that both the co-ordinates of $P_2$ are positive. 

Thus, starting from a known point $P_1$ on the curve \eqref{cubic2} given by $(x_1,\,y_1)$ with $2^{2/3} < x_1 <  2$ and $y_1 >0$, we have generated a new point $P_2$ on \eqref{cubic2} given by $(x_2,\,y_2)$ with $2^{2/3} <  x_2 <  2$ and $y_2 > 0$. We now start with the point $P_2$ and repeat the process described above to get a new point $P_3$ whose co-ordinates $(x_3,\,y_3)$ satisfy the relations $2^{2/3} < x_3 <  2$ and $y_3 > 0$. In fact, we can  execute  this process repeatedly to get an infinite  sequence of points  $P_1,\;P_2,\;P_3,\;\ldots,P_n,\ldots,$ on the curve \eqref{cubic2} given by $(x_i,\,y_i),\;i=1,\,2,\,3,\ldots,n,\,\ldots$  such that for each $i$, we have  $2^{2/3} <  x_i <  2$ and $y_i > 0$. We will now show that with a suitably chosen initial point $P_1$, all the points of this sequence are distinct.

If we take the point $P_1$ as $n_1P$ where $ 2 < |n_1|$, then the point $Q =2P_1-P = (2n_1-1)P$, and so we get $P_2=n_2P$ where $n_2$ is either $2n_1-1$ or $-2n_1+1$. Thus, $ 2 < |n_1| < |n_2|$, and, by induction, it follows that the sequence of points $P_1,\;P_2,\;P_3,\;\ldots,P_n,\ldots,$ is given by $n_1P,\;n_2P,\;n_3P,\ldots,$ where the absolute values of $n_i$ form a strictly increasing sequence of positive integers. Since the point $P$ is a point of infinite order, we are assured of an infinite sequence of distinct rational points $(x_i,\,y_i)$ on the curve \eqref{cubic2} with   $2^{2/3} <  x_i <  2$ and $y_i > 0$. Finally, we use the relations \eqref{birat1} to obtain a corresponding sequence of infinitely many points on the curve \eqref{cubic1} whose co-ordinates are positive rational numbers.

As an example, if we take the point $P_1$ as the point $A$ given by $(785/484,$  $5497/10648)$ which is $-4P$, we obtain an infinite sequence of points  $P_1,\;P_2,$ $P_3,\;\ldots,P_n,\ldots,$ with the specified properties.  The point $P_2$ of this sequence is  given by
\begin{multline}
(8152570498330546/4944742493612769, \\
241351355149002573947470/347708669978634678361647).
\end{multline} 
The next point of the sequence, that is $P_3$, has co-ordinates involving integers with 85 digits and is therefore omitted. The solution of \eqref{eqk4} corresponding to the point $P_2$ is given by,
\begin{equation}
\begin{aligned}
a& = 234192173776567982667691/286639743984973696444599,\\
b& = 113516496202066695693956/286639743984973696444599.
\end{aligned}
\end{equation}

\subsection{} Finally we note that when $k \geq 5$, Eq.~\eqref{eqk} represents a curve of genus $> 1$, and by Falting's theorem, it can have only finitely many rational solutions. Limited trials conducted on Eq.~\eqref{eqk} with $k=5$ yielded no nontrivial solutions. It appears that Eq.~\eqref{eqk} has  no nontrivial rational solutions when $k \geq 5$.

\section{Numerical Curios Involving Radicals}

We now show how to obtain infinitely many numerical curiosities involving radicals and rational numbers of the type \eqref{ex1}. It is  easily seen that 
any solution of Eq.~\eqref{eqk} leads to
\begin{equation}
\sqrt[k]{b^k}+d=
\sqrt[k]{b^k+d}
\quad{\hbox{and}}\quad
\sqrt[k]{a^k}-d=
\sqrt[k]{a^k-d}, \label{curiogen}
\end{equation}
where $d=a^k-b^k=a-b$.

When $k=2$, the simple solution $a=2/3,\;b=1/3$ of Eq.~\eqref{eqk} yields the examples,
\[
\sqrt{\frac{1}{9}}+{\frac{1}{3}}
=\sqrt{\frac{1}{9}+{\frac{1}{3}}} \ \phantom{.}, \quad \mbox{\rm and} \quad  \sqrt{\frac{4}{9}}-{\frac{1}{3}}
=\sqrt{\frac{4}{9}-{\frac{1}{3}}}\ \phantom{.}, \]
while the solution \eqref{solDiophantus} of Eq.~\eqref{eqk} when $k=3$ leads to the identities,
\[
\sqrt[\scriptstyle3]{\frac{343}{2197}}+{\frac{1}{13}}
=\sqrt[\scriptstyle3]{\frac{343}{2197}+{\frac{1}{13}}}\ \phantom{.},
\]
and 
\[
\sqrt[\scriptstyle3]{\frac{512}{2197}}-{\frac{1}{13}}
=\sqrt[\scriptstyle3]{\frac{512}{2197}-{\frac{1}{13}}} \ \phantom{.} 
.
\]

The solution \eqref{sol1eqk4} of Eq.~\eqref{eqk4} given by Fermat yields the identity \eqref{ex1} given at the beginning of this paper. The infinitely many solutions in positive rational numbers of Eq.~\eqref{eqk} when $k=2,\,3,$ or 4 yield infinitely many similar identities.

\section{Geometric Series with an interesting property}

Solutions of the  diophantine equation \eqref{eqk} also enable us  to construct  examples of  convergent geometric series $\sum_{n=1}^\infty a_n$ with positive rational terms and possessing the interesting property  described by Eq.~\eqref{series}.

In view of the formula for the sum of infinite geometric series, Eq.~\eqref{series}
can be presented as
\begin{equation}
\frac{a_1}{1-r}=
\frac{a_1^k}{1-r^k}
,
\label{series-eq}
\end{equation}
where $r$ is the common ratio of the geometric series $\sum_{n=1}^\infty a_n$.
Now Eq.~\eqref{series} can be easily rewritten as
$$
\frac{1}{a_1}-
\frac{q}{a_1}=
\left(\frac{1}{a_1}\right)^k-
\left(\frac{q}{a_1}\right)^k
,$$
which is equivalent to the equation (\ref{eqk}) by substituting
$a=1/a_1$ and $b=r/a_1$ or equivalently,
$a_1=1/a$ and $r=b/a$. To ensure that the two geometric series are convergent,  we choose $a,\,b$ such that $b < a$, and hence $r < 1$. 

The infinitely many solutions of \eqref{eqk} with $k=2,\,3$ or $4$ yield infinitely many examples of geometric series satisfying Eq.~\eqref{series}. As a numerical example, the solution \eqref{solDiophantus} of Eq.~\eqref{eqk} when $k=3$ leads to the two geometric series,
\begin{equation}
\left(\frac{\displaystyle 13}{\displaystyle 8}\right)+\left(\frac{\displaystyle 13}{\displaystyle 8}\right).\left(\frac{\displaystyle 7}{\displaystyle 8}\right)+\left(\frac{\displaystyle 13}{\displaystyle 8}\right).\left(\frac{\displaystyle 7}{\displaystyle 8}\right)^2+\left(\frac{\displaystyle 13}{\displaystyle 8}\right).\left(\frac{\displaystyle 7}{\displaystyle 8}\right)^3+\cdots,
\end{equation}
and
\begin{equation}
\left(\frac{\displaystyle 13}{\displaystyle 8}\right)^3+\left(\frac{\displaystyle 13}{\displaystyle 8}\right)^3.\left(\frac{\displaystyle 7}{\displaystyle 8}\right)^3+\left(\frac{\displaystyle 13}{\displaystyle 8}\right)^3.\left(\frac{\displaystyle 7}{\displaystyle 8}\right)^6+\left(\frac{\displaystyle 13}{\displaystyle 8}\right)^3.\left(\frac{\displaystyle 7}{\displaystyle 8}\right)^9+\cdots,
\end{equation}
both of which have the same common sum, namely 13. 

As a second example, 
 the solution $a=26793/34540$, $b=15799/34540$ of Eq.~\eqref{eqk4} leads to the geometric series defined by
$$a_n=\frac{34540\cdot15799^{n-1}}{26793^n}
$$
satisfying the condition
$
\sum_{n=1}^\infty a_n=
\sum_{n=1}^\infty a_n^4=17270/5497
$.

\noindent Postal Address 1: Ajai Choudhry, \newline \hspace{1.05 in}
13/4 A Clay Square,
\newline \hspace{1.05 in} Lucknow - 226001, INDIA.
\newline \noindent  E-mail: ajaic203@yahoo.com

\noindent Postal Address 2: Jaros\l aw Wr\'oblewski, \newline \hspace{1.05 in}
Mathematical Institute, Wroc\l aw University,
\newline \hspace{1.05 in} pl. Grunwaldzki 2/4,
\newline \hspace{1.05 in} 50-384 Wroc\l aw, POLAND.
\newline \noindent  E-mail: jwr@math.uni.wroc.pl

\end{document}